\documentclass[]{article}

\usepackage{amssymb, epsf, lscape, color, colortbl, subfigure, epsfig, setspace,graphicx,slashbox}
\usepackage{amsmath, amsthm, fancybox, float, amsfonts}
\usepackage[english,]{babel}
\usepackage{makeidx}
\newcommand{\lct}{allocate }
\newcommand{\lcn}{allocation }

\newcommand{\dst}{distribution }

\newcommand{\dstr}{destruction }
\newcommand{\dsts}{distributions }

\newcommand{\flw}{following }
\newcommand{\fls}{follows}

\newcommand{\crsp}{corresponding}

\newcommand{\ndp}{independent}
\newcommand{\ndpl}{independently }

\newcommand{\prb}{probability }
\newcommand{\prbe}{probability}

\newcommand{\prbs}{probabilities }

\newcommand{\stg}{strategy }
\newcommand{\stge}{strategy}
\newcommand{\stgs}{strategies }
\newcommand{\stgse}{strategies}

\newcommand{\ds}{\displaystyle}

\makeatother
\begin{document}

\author{Isaac M. Sonin\\
University of NC at Charlotte}
\title{Bayesian Game of Locks, Bombs and Testing}
\date{\today }
\maketitle

\thispagestyle{empty}\bigskip \bigskip




\section{General Bayesian Game of Locks, Bombs and Testing}\label{Lbt}
     \large {
In this note, we describe General Bayesian Game of Locks, Bombs and Testing, briefly LBT game,
and the main points of the two incoming papers, \cite{ls19} and \cite{ss19}, where some of important special cases were solved. The inspiration for this model was the paper by K. Sonin, J. Wilson, A. Wright: Rebel Capacity, Intelligence Gathering, and Combat Tactics, (manuscript, \cite{sww18}).\par
In a classical  Blotto game, two players distribute limited resources between different sites (battlefields) with the goal to win more sites, winning a site if you have more resources on this site than your opponent. There is substantial literature on this topic, with classic paper \cite{shub81} and more recent publications, such as \cite{rob06}, where a complete solution of the “continuous” version was given, as well as \cite{har15}, where some interesting extensions are discussed. In a comprehensive and detailed survey \cite{hoh16} dedicated to Search Games, the Blotto game is classified as an attack-defense game. There are even more papers dedicated to these games and, as in all of Operations Research all classifications have many overlapping parts. As an example of an important paper of an attack-defense game we mention \cite{pow07}. Our model shares some similarities with the model in this paper and the other papers about the Blotto game, but also features some important distinctions.\par
In our model, as in most attack-defense games, the two players are not symmetrical. We call one of them Defender (DF) and the other, Attacker (AT). There are $n$ ``sites"\ (battlefields, boxes, cells, targets, time slots, etc.) with a \emph{vector of values} (costs) $c=(c_i,i=1,...,n)$. AT is trying to destroy these sites by placing ``bombs"\ that can result in explosions and destruction, and DF is trying to protect these sites by distributing ``locks"\ among them. One or more bombs can be placed into the same site. A lock is a protection device which, placed in a site, prevents its \dstr with any number of bombs in it. Thus, locks and bombs are just the names of discrete units of resources of protection and destruction. The number of locks $k, k<n$, can be fixed (paper \cite{ss19}) or it can be a random variable, obtained, e.g., if a lock appears in site $i$ with \prb $\lambda_i$ \ndpl of other sites. The latter case is the subject of paper \cite{ls19}. Here we discuss only the case of $k$ locks. AT has $m, m=1,2,...,$ bombs to \lct among $n$ boxes. A box is destroyed if at least one explosion occurs, and the explosions of different bombs in the same site or in different sites are \ndp. We denote $p$ the \prb of explosion.\par
The important feature of our model is that AT can and will \emph{test} every site, trying to find sites without locks. This testing is not perfect: a test of site $i$ may have a positive result, $S_i=1,$ even if there is no lock at the site, $T_i=0$,  and negative, $S_i=0$, even if there is a lock, $T_i=1$. The probabilities of correct identification of both types, in a statistical language the \emph{sensitivity} and \emph{specificity}, $P(S_i=1|T_i=1)=a_i$ and $P(S_i=0|T_i=0)=b_i$, are known to both players. The result of testing is a vector of signals $s=(s_1,...,s_n)$, where each $s_i=0,1$, known to AT.\par
When the number of locks $k, k<n$, is fixed, then the \emph{strategy of the defender} (DF) is a \prb \dst
$b(\gamma )$ on a set of all possible positions of locks $\{\gamma\}$, where $\gamma =(i_{1},i_{2},...,i_{k})$ with $1\leq i_{1}<...<i_{k}\leq n$. In our Bayesian setting we assume that this prior \dst $b(\gamma )$ \emph{is known} to AT, although of course the actual positions of the locks are not. After the locks are allocated, AT receives signal $s$ and, having $m$ bombs, distributes them among $n$ sites deterministically or using some randomization, trying to maximize
\emph{the expected sum of values of all destroyed sites}. WLOG, we can assume that this \dst is deterministic and an \emph{optimal \stg of AT} $\pi(m | b(\gamma ))$, with respect to a \stg of DF $b(\gamma )$, is a collection of her optimal responses $u(s|m, b(\gamma ))\equiv u(s)=(u_1(s),...,u_n(s))$ to each signal $s$, where $u_i(s)$ is the number of bombs placed into site $i, i=1,...,n$. Using prior \dst $b(\gamma )$, the \prbs of signals $p(s)\equiv p(s|b(\gamma ))$, the aposterior \dst of the positions of locks $b(\gamma|s)$, the optimal \stg $\pi(m |b(\gamma ))$ and the \crsp \ total expected damage (loss), $L(b(\gamma )$, can be calculated. The goal of DF is to select a prior \dst of locks $b_*(\gamma )$ to minimize this loss. Then the pair $(b_*(\gamma ), \pi_*)$, where $\pi_*$ is an optimal response of AT to \stg $b_*(\gamma )$, forms a classical Nash equilibrium (NE) point. The corresponding value of the game is $v_*=L(b_*(\gamma ), \pi_*))$. Though $b_*(\gamma )$ are not unique, they all have common properties that result in a unique (up to some randomization) AT \stg $\pi_*$, and thus a specific value of $v_*$. We denote $G\equiv G(n,k,m|a,b,c)$ this general Locks, Bombs and Testing (LBT) game (model), where $n$ is the number of sites, $k$ is the number of locks, and $n$-dimensional vectors $a=(a_i),b=(b_i)$ and $c=(c_i), i=1,2,...,n$, represent the sensitivities and the specificities of testing, and the values of all sites. \par
Note that \emph{sites, locks, bombs  and testing} in this model are rather abstract terms and may have very
different interpretation beyond our initial exposition of the DF and AT defence-attack game. For example, the American presidential elections with the distribution of efforts [resources] between ground
battle states,  with polls serving as testing tools, fit this model. In biology and medicine, sites may represent parts
of the body (organs) and bombs units of treatment (chemo, radiation) delivered to metastatic points and subjected to limitation for a total dose (amount). In Search theory, locks are hidden objects and bombs are resources spent on detection of the objects. In repair/maintenance  models one may consider locks as hidden defective blocks and bombs as some maintenance/repair units. The number of bombs in an $i$-th site (block) defines the \prb of repair of this site. Thus, the described game is very flexible and allows generalizations in many directions: continuous resources, diverse locks and bombs, fake locks, etc. For example, we can assume that locks can not protect sites from destruction but can distort signals, and then such model will have common features with modern Bayesian persuasion models. \par
\vspace{.2cm}
There are many possible modifications of the LBT game. We mention only some of them. We assumed that the parameters $a_i,b_i, i=1,...,n$ are known to DF. Without this assumption the solution should be quite different.
We assumed that the prior \dst of locks is known to AT. This assumption can be somewhat relaxed but the discussion of this subject leads us to the deeper levels of Game Theory.\par
Finally, we would like to mention that the \emph{dynamic version }of the GLBT model will have common features with the well known model in Applied Probability - the Multi Armed Bandit problems. This model was studied in many papers and a few books - D. Berry and B. Fristed (1985), E. Presman and I. Sonin (1987, 1990), J. Gittins (1989), J. Gittins, K. Glazebrook and R. Weber (2011), and the current internet version by T. Lattimore and C. Szepesv´ari (2019).   \par
The full solution of the LBT game is a difficult tusk. In the next section, we describe the full solution of the this game for the so called symmetrical case, in Section 3 we present some partial results for the general case, and in Section 4 we present the full solution for the case $G(2,1,1$.
\par
\section{Symmetrical LBT model}
\subsection{Setup}
The general $G$ game and the straightforward approach to solving it, described above, have two basic drawbacks. First, the set of possible positions for locks, i.e., the set of subsets of an $n$ element set, has order $2^n$, and so does the set of potential signals. As a result, the calculations of posterior \dsts $b(\gamma|s)$ and their marginal \dsts $\alpha_i(s)=P(T_i=0|s)$ which play a crucial role in the description of optimal \stgse, become cumbersome for large $n$. The second problem is that the knowledge of detailed information about the values of $c_i, a_i$, and $b_i$ in many cases is unrealistic. As a result, the main focus in papers \cite{ls19} and \cite{ss19} was on the analysis of a simpler  \emph{symmetrical LBT (S-LBT) model}, where all sites have identical values $c_i=1$, and all $a_i=a, b_i=b$. \par
Note that both general and symmetrical models have one special and important feature. Its analysis and the solutions consist of two parts. In the first part, we consider a \emph{statistical problem}: in general model, we need to analyze the aposterior \dst of locks (ADL) $b(\gamma|s)$, and in symmetrical model, it is sufficient to analyze only the two critical ratios defined in this section. In the second part, we have to solve an \emph{optimization problem }of bombs allocation using signal $s=(s_1,...,s_n)$, the aposterior \dst of locks, and the explosion function $p(u)$, i.e., the \prb of at least one explosion of $u$ bombs. \par
In both general and symmetrical models, there are also two possible assumptions about the number of locks $k$: problem $A(n,k)$, where the number of locks $k$ is fixed, and problem $B(n,\lambda)$, where each of $n$ sites has a lock with \prb $\lambda$ \ndpl of other sites. In S-LBT model, DF, having no information about AT, uses the uniform \dst of locks, and we have a one-sided problem to find an optimal \stg for AT. In problem $A(n,k)$, this uniform prior \dst in statistical physics is called Fermi-Dirac statistics and any combination of $k$ protected boxes has the same \prb $\ds 1/\binom{n}{k}$. It is easy to show that the \prb of protection for a particular box is $\ds t=\frac kn$. In problem $B(n,\lambda)$, the number of locks has a binomial \dst with parameters $n,\lambda$. \par
In both S-LBT models, in the first, statistical part of the solution, given signal $s=(s_1,...,s_n)$, instead of using the aposterior \dst of locks, it is sufficient to analyze only two critical ratios: ratio $r$ for problem $B(n,\lambda)$, and $r(x)$ for problem $A(n,k)$:\par
\vspace{-.2cm}
\begin{eqnarray}
r=\frac{P(T_i=0|S_i=0)}{P(T_i=0|S_i=1)}\equiv\frac{p^{-}}{p^{+}},\ \ \ \ \ r(x)=\frac{P(T_i=0|S_i=0,N=x)}{P(T_i=0|S_i=1,N=x)}\equiv\frac{p^{-}(x)}{p^{+}(x)}. \label{rxr}
\end{eqnarray}
As a result, the optimal \stg in both problems will have a much simpler structure than in the general case, symmetrical with respect to all sites with negative (minus) signals, and correspondingly for sites with positive (plus) signals.
In this statistical part, the explosion function $p(u)$ and the \stg of bombs allocation do not participate at all.\par
Though the solutions of both problems have certain similarities, some of their features are very distinct. For example,
an interesting and even counterintuitive property is that in problem $A(n,k)$, the function $r(x)$ and
therefore the optimal \stg and the value function, depend on $a$ and $b$ only through the value $\ds q=\frac{a}{1-a}\frac{b}{1-b}$, a combined characteristic of the quality of testing. In problem $B(n,\lambda)$, this property is not true with respect to the value $r$. The other important distinction between these two problems is that in the former problem the minuses and pluses in different boxes are not independent, but in the latter problem they are. \par
In both problems, DF and AT know all relevant parameters and prior \dst $b(\gamma)$, but only AT knows the signal. The reader is advised to consider the outcomes as a result of a two-stage random experiment. In problem $B(n,\lambda)$, at the first stage, a lock appears in every site independently of other sites with \prb $\lambda$. In problem $A(n,k)$, $k$ locks are distributed at random between $n$ sites, where $k$ is a fixed number. After that, at the second stage, in both problems, all $n$ sites are tested. The results of the tests are independent from site to site and depend only on whether a lock is in the site or not.
Formally, we consider random variables $T_i, S_i, C_i, i=1,2,...,n$, taking two values, $0$ and $1$; $T_i=1, (0)$ means
that the $i$-th site has a lock, (no lock); $S_i=1$ or $+$, ($0$ or $-$), means that the test of the $i$-th site is positive, (negative), $C_i=1$, ($0$) means that the $i$-th site is destroyed, (not destroyed). The absence of subindex $i$ means that we discuss an arbitrary site, e.g., $P(T=0|S=0)$, etc.
The outcomes of the two-stage random experiment described above, i.e., the elements of the \prb space, are pairs $(\gamma, s)$, where $\gamma_k\equiv \gamma=(i_1,i_2,...,i_k)$ with $1\leq i_1<...<i_k\leq n$, is a (vector) \emph{position of $k$ locks} and $s=(s_1,...,s_n)$ is a (vector) \emph{signal}. The \prb of each outcome is $P(\gamma, s)=b(\gamma)P(s|\gamma)$, where $b(\gamma)$ is the prior \dst of locks and $P(s|\gamma)=P(S_1=s_1,...,S_n=s_n|\gamma)$. \par
\vspace{-.2cm}
\subsection{Optimal Strategies and Values in a Symmetrical LBT Model}\label{Osv}
To solve LBT models means to find an optimal \stg and the corresponding value function, i.e., the value of the functional under the optimal \stge. A \stg $\pi\equiv \pi(\cdot |m,s)=(u_1, u_2, \ldots, u_{n}|m,s)$,  $\sum_{j=1}^{n}u(j)=m$, is an allocation of $m$ bombs in $n$ sites, given signal $s$, defined for all $m$ and $s$. We will skip sometimes to indicate dependence on $m$.\par
The symmetry in the S-LBT model implies two useful formulas
\vspace{-.3cm}
\begin{eqnarray}
P(s_1,...,s_n)=P(N=x)/\binom{n}{x}, \ \ P(T_i=0|s_1,...,s_n)=P(T_i=0|s_i,N=x). \label{8}
\end{eqnarray}
The first equality is almost obvious, since the symmetry of signals and bombs \lcn implies that $P(s|N=x)=c(x)$, i.e., all signals with the same number of minuses $x$ have the same \prbe. Then, since the number of such signals is $\binom{n}{x}$ and $\sum_sP(s|N=x)=P(N=x)$, we obtain that $P(s)$ is given by the first formula in (\ref{8}). The second formula in (\ref{8}) can be obtained from the first one used for the cases of $n$ and $n-1$ sites.
Let us denote $B^{-}(s)$ and $B^{+}(s)$ the sets of minuses and pluses, given signal $s$. Then, using the second equality in formula (\ref{8}), and the definitions of $p^{-}(x),p^{+}(x),r(x)$ given in formula (\ref{rxr}), we obtain that, given signal $s$ with $N(s)=x$, and a \stg $\pi=(u_1,...,u_n|s)$ for $m$ bombs, the value of a \stg $\pi$ is
\vspace{-.2cm}
\begin{eqnarray}
w^\pi(s)\equiv w^\pi(s,m)&=&\sum_{i=1}^{n}P(T_i=0|s)p(u_i)=\sum_{i=1}^{n}P(T_i=0|s_i,x)p(u_i)= \notag \\
&=&p^{+}(x)[r(x)\sum_{i\in B^{-}(s)}p(u_i)+\sum_{i\in B^{+}(s)}p(u_i)].  \label{ws}
\end{eqnarray}
Let $U^{-}\equiv U^{-}(\pi|s)=\{u_j, j\in B^{-}(s)\}$ and $U^{+}\equiv U^{+}(\pi|s)=\{u_j\in B^{+}(s)\}$ be two possible sets of the values of $u_j$ in minus and plus sites. Formula (\ref{ws}) immediately implies that in the S-LBT model all \stgs obtained by permutations of sets $(U^{-}, U^{+})$ among the \crsp \ sites have the same value, denoted as $w^\pi(x)\equiv w^\pi(x,m)$.
In other words, a \stg in the S-LBT model can be understood as an \lcn of bombs in minus and plus sites without taking into account the particular numbers of sites.\par
We also denote $v(x,m)\equiv v(x)=\sup _{\pi}w^\pi(x,m)$, the \emph{value function} over all such strategies, given $m$ and $x$, and $v(m)$, the value function over all \stgs and all possible values of $x$, i.e., $v(m)=\sum_xP(N=x)v(x,m)$.\par
Formula (\ref{ws}) gives a hint that the proportion of the number of bombs placed into a minus site to the number of bombs placed into a plus site is defined by ratio $r(x)=p^{-}(x)/p^{+}(x)$, defined in formula (\ref{rxr}).
Under a ``normal"\ situation, when $a$ and $b$ are such that tests provide ``truthful"\ information, $r>1$ and $r(x)>1$ for all $x$, and therefore, if there is only one bomb, it should be placed into a minus site. We show later that ``normal"\  means that $a+b>1$. When the number of bombs $m$ exceeds $x$, optimal \stgs can be expressed through $r$, $r(x)$, and other parameters. In a sense, the values $N=x, r, r(x)$ play the role of sufficient statistics in the optimization problem. \par
Given $N=x, 0\leq x\leq n$, the \lcn of bombs depends on the number $m$ of the bombs available. The optimal \stg has the \flw structure. Initially, all bombs are allocated one by one into each of the $x$ minus sites until the level $d(x)\geq 1$ is reached in each of them or the bombs are exhausted. Afterwards, the bombs are added one by one to plus sites to reach level one in each of them. Afterwards, the bombs are added one by one to minus sites to reach level $d(x)+1$ in each of them, then back to plus sites to reach level two, etc. Stop this ``fill and switch"\ process when run out of bombs. We called such \stg $d(x)$ ``uniformly as possible"\ $d(x)$-(UAP) \stge, and this process implies that either all minus sites have the same level, or all plus sites have the same level, or both, and if all plus sites are empty, then the level in minus sites $l^{-}\leq d(x)$, and if some plus sites have bombs, then the difference between the levels equals either $d(x)$ or $d(x)-1$. Let us denote $m^{-}, m^{+}$ the numbers of bombs in the minus and plus sites. In an example with $n=5, x=3, m=10$ and $d(x)=2$, $d(x)$-UAP implies that $m^{-}=8, m^{+}=2$, two minus sites have three bombs each, one minus site has two bombs, and two plus sites have one bomb each. If $x=0$ or $n$, then all sites are filled sequentially one by one and we call this \lcn just UAP. \par
Formally, given $m$ and $x$, $0<x<n$,  $d(x)$-\stg defines a unique \lcn of bombs, given by tuple $(l^{-}, e^{-}, l^{+}, e^{+})$, where $l^{-}, (l^{+}$) is the number of ``complete layers"\ of bombs in the minus (plus) sites, and $e^{-}, (e^{+})$ is the number of extra bombs in the ``incomplete layers"\ in the minus (plus) sites. All these terms depend on $m,x$ and $d(x)$, but we do not indicate this explicitly. Using shorthand notation $l^{-}=i, e^{-}=e,\ \ l^{+}=j, e^{+}=e'$, we have
$m=m^{-}+m^{+},\ m^{-}=i*x+e,\ m^{+}=j*(n-x)+e',  0\leq e<x, 0\leq e'<n-x$ and $e*e'=0$.
Thus, if $m^{+}=0$, then $m=m^{-}\leq d(x)$;\ if $e'>0$, then $e=0$ and $i-j=d(x)$; if $e'=0, j>0$, then either $i-j=d(x)-1, e\geq 0$ or $i-j=d(x), e=0$.\par
\vspace{.2cm}  
 \textbf{Theorem 1.} \emph{Let, given signal $s$, the total number of minuses $N=x, 0\leq x\leq n$. Then} \par
a) \emph{if $x=0$ or $n$, then the optimal \stg is to distribute all bombs between the sites UAP and the value function $v(m|0)=v(m|n)$ for $m=n*i+e, \ i=0,1,...,\ 0\leq e<n$, is given by formula}
\vspace{-.2cm}
\begin{eqnarray}
v(m|0)&=&v(m|n)=\frac{n-k}{n}[ep(i+1)+(n-e)p(i)]; \ \  \label{vm0}
\end{eqnarray}
b) \emph{if $0<x<n$, then the optimal \stg is to distribute all bombs between the minus and the plus sites $d(x)$-UAP,
where $d(x)$ is defined by formula  }
\vspace{-.5cm}
\begin{eqnarray}
d_{}(x)=min(i\geq 1: r(x)q^i<1),  \label{dx}
\end{eqnarray}
\emph{$q=1-p$, and $r(x)$ is defined by formula} (\ref{rxr}).\par
\emph{The value function $v(x,m)$ for $m=m^{-}+m^{+}=i*x+e+j*(n-x)+e'$, where the tuple $(i,e,j,e')$ is (uniquely) defined by value $x$ and a $d(x)$-UAP \stge, is given by formula}
\vspace{-.2cm}
\begin{eqnarray}
v(x,m)=p^{+}(x)[r(x)(ep(i+1)+(x-e)p(i))+e'p(j+1)+(n-x-e')p(j)].\label{vmx}
\end{eqnarray}
c) \emph{The value function $v(m), m=1,2,...$, is given by formula}
\vspace{-.2cm}
\begin{eqnarray}
v(m)=\sum_{x=0}^{n}P(N=x)v(x,m). \label{vm}
\end{eqnarray}
\textbf{Remark 1.} \emph{By definition of $d(x)$, we have $r(x)q^{d(x)-1}\geq 1$. If $r(x)q^{d(x)-1}>1$, then $d(x)$-UAP \stg is the unique optimal \stge. If there is the equality, then there are other optimal \stgs with the \lcn of bombs obtained as follows. When all minus sites are filled with $d(x)-1$ full layers, the next bomb, if available, can be placed either in a minus site or in a plus site. The incremental utility will be the same. And so on with other extra bombs. As a result, the difference between the full layers in the minus, and the plus sites can be either $d(x)$ or $d(x)-1$.} \par
\vspace{.1cm}
\section{Partial Results for the General LBT Game}
  Before formulating the most general results, we consider some relatively simple situations. If testing is \emph{perfect}, i.e., when  $a_i=b_i=1, i=1,...,n$, the optimal \stge for AF is clear: having $k$ locks, she has to defend the $k$ most valuable sites, placing into them all $k$ locks. With only one available bomb, $m=1$, AT will place it into the next valuable site, and if $m>1$, she should solve the problem of discrete optimization, each time placing the next available bomb into the site with the maximal marginal utility.\par
The other extreme situation is when testing is \emph{non-informative}, i.e., when  $a_i=b_i=1/2, i=1,...,n$. Then the aposterior \dst $b(\gamma|s)$ coincides with the prior \dst $b(\gamma )$ for all signals $s$. Let us say $n=2, k=1, m=1, c_1=c, c_2=1, p=1$. At first glance, it seems that if $c$ is much larger than 1, DF should place a unique lock into the most valuable site, and then her loss is 1. But simple calculations show that the optimal \dst of locks is $(\frac{c}{c+1}, \frac{1}{c+1})$ and $v=\frac{c}{c+1}<1$. AT can place her unique bomb into any site, or place it at random. \par
Similarly, for the game $G(n=3,k=1, m=1)$ with vector of values $c=(4, 3, 2)$, it is easy to obtain that the optimal \dst of a unique lock is given by $b(\gamma)=(\frac{7}{13}, \frac{5}{13}, \frac{1}{13})$ and $v_*=\frac{24}{13}\approx 1.85$. But if the vector of values is $c=(4, 3, 1)$, then $v_*=\frac{12}{7}\approx 1.71$ and $b(\gamma)=(\frac{4}{7}, \frac{3}{7}, 0)$. Thus, under the optimal \stg of DF, site 3 is not protected at all. \par
Given prior \dst $b(\gamma)$, let us introduce marginal \prbs $\beta_i=P(T_i=1)$ and their complimentary \prbs $\alpha_i=P(T_i=0)$. It is easy to show that $\sum_{i=1}^n\beta_i=k$ and that $\sum_{i=1}^n\alpha_i=n-k$, i.e., equal 2 in our examples. Then in the first example above, we have vector $\alpha=(\frac{6}{13}, \frac{8}{13}, \frac{12}{13})$ and thus $c_i\ast\alpha_i=v_*=\frac{24}{13}, i=1,2,3$. In the second example, we have $\alpha=(\frac{3}{7}, \frac{4}{7}, 1)$ and thus $c_i\ast\alpha_i=v_*=\frac{12}{7}, i=1,2$ and $v_*>c_3=1$. In both examples, the optimal \stg of DF is unique, but for larger values of $n$ this is not the case, though all optimal $b_*(\gamma)$ have the same \prbs $\alpha_i=P(T_i=0)$. A similar situation holds in the general case. We present here the theorem for the game $G(n, k, m=1)$. We present the full version of this theorem for any $m$ and its proof in the extended paper.\par
\vspace{.2cm}
\textbf{Theorem 2. (Non-Informative Case)}. For $m=1, c_1\geq c_2\geq ...\geq c_n$, \par
a) the class of optimal \stgs $b_*(\gamma)$ has the following structure:  there is a $k_*=k_*(c), k<k_*\leq n$ and a constant $v_*=v(c)$, such that: $c_i\alpha_i=v_*$ for $1\leq i \leq k_*$, and $v_*>c_i, \alpha_i=1$  for $k_*<i\leq n$;\par
b) the optimal \stg for AT is to place a bomb at random between the sites with numbers $1,2,...,k_*$;\par
c) the value of the game is $v_*=(k_*-k)/C_{k_*}$, where $k_*=max\{j\geq k: c_j>(j-k)/C_j\}$ and $C_j=\sum_{i=1}^{j}1/c_i, j=1,2,...,n$.
\par
\vspace{.2cm}
In other words, if $k_*<n$, then the sites with numbers greater than  $k_*$ should not be protected at all, and the \dst of $k$ locks in the first $k_*$ sites should make all sites \emph{equally desirable for attack}.  For a Bayesian zero-sum game, this equality is a well-known fact called \emph{the Principle of Indifference}, (see e.g. \cite{mson91} or \cite{ferg14}.) Citing the latter text, ``Player 1 searches for a strategy that makes Player 2 indifferent
as to which of the (good) pure strategies to use". The similar statement for Player 2 in our case has a different character.
In layman terms: if the strength of a chain is defined by the strength of the weakest link, and the resources to make links strong are limited, then make all links of equal strength.  This is a special case of a more general ``Chain-Link Optimization Principle" that can be applied to many other optimizations problems, (\cite{son19}). \par
Note also that for \emph{any vector of values} $c$, the number $k<k_*$, the optimal \stg $b_*(\gamma)$ is always randomized and value $v_*>c_{k_*+1}$. As a result, AT will allocate her $m$ bombs among the first $k_*$ sites if $m \leq k_*$. \par
An interesting situation appears when there is a mixture of full informative and non informative sites, i.e., when $a_i=b_i=1$ for some sites and  $a_i=b_i=1=1/2$ for all others. Then a theorem similar to Theorem 1 can be proved with an explicit description of optimal the \stgs for both players. \par
Of course, the main interest in the Bayesian LBT game problem is the case of \emph{imperfect but informative} testing. This means that $1/2<a_i, b_i<1, i=1,...,n$.
To obtain the description of NE points, we have to solve three problems. \par
\textbf{The first problem} is, given a \stg of DF $b(\gamma )$, to describe the optimal \stg (response) of AT $\pi(m | b(\gamma))$, i.e., to describe the optimal \lcn of bombs $u(s|m)$ given signal $s$ and $m$ available bombs. The full answer to this problem is given by a recursive procedure $S$ described in \cite{ss19}.
The expected value of damage (loss) for the pair of strategies $(b(\gamma),\pi(m | b(\gamma))$ can also be obtained.\par
\textbf{The second}, more difficult, problem is to find the optimal \stg or \stgs $b_*(\gamma)$ of DF, minimizing this loss. 
The heuristic meaning of the \crsp \ Theorem 3 is similar to the meaning of Theorem 1: these strategies have to make the potential expected losses in the sites, that are worth protecting, equal when AT applies her optimal response to $b_*(\gamma)$. \par
\emph{Given a \stg of DF} $b(\gamma )$, let us denote $L_i(s|m)$ the expected loss in site $i$ when AT applies her optimal \stg given signal $s$, $p(s)$ the \prb of signal $s$ and $L_i(m)=\sum_{s}p(s)L_i(s|m)$ is the \crsp \ expected loss. \par
\vspace{.2cm}
\textbf{Theorem 3. (Informative Case).} Given $m=1,2,...$, the class of optimal \stgs $b_*(\gamma|m)$ can be obtained using the Principle of Indifference, that takes the \flw specific form:
the values of $L_i(m)=v_*$ must be equal for all $i=1,2,...,k_*(m)$, where $k_*(m), k<k_*(m)\leq n$, is the number of sites worth protecting, and $v_*$ is the value of the game.\par
\vspace{.2cm}
The difficulty in applying Theorem 3 lies in the fact that in the informative case the optimal response depends on signal $s$, and the calculation of $L_i(m)$ is nontrivial. In the next section, we provide an example of an application of Theorem 3 to find the optimal \stge.\par
\textbf{The third problem} to be solved is to obtain the full description of all NE points, i.e., to describe all $b_*(\gamma)$ delivering the equality of $L_i(m)$ in Theorem 3.
\par
\vspace{.2cm}
\textbf{Remark 1}. The description of $b_*(\gamma )$ is based on the \flw interesting property of the $G$ game: to obtain the optimal response of AT given any $b(\gamma )$ and signal $s$, i.e., to use procedure R, AT needs to know only the\emph{ marginal \prbs} $\alpha_i(s)$ for all $s$, but DF, trying to obtain $b_*(\gamma )$, needs to know $v(m)$ and the total expected loss, so she needs to know $p(s)$, which are based on the \emph{whole \dst }$b(\gamma )$. There is a simple example that shows that two distinct $b(\gamma )$ can have the same  \prbs $\alpha_i(s)$. Thus, one of the side problems is to obtain the description of all $b(\gamma )$ having the same \prbs $\alpha_i(s)$. \par
\vspace{.2cm}
\textbf{Remark 2}. The statements and interpretation of Theorems 1 and 2 can be expressed also using the concepts of \emph{information} and \emph{entropy}. Loosely speaking, the optimal DF \stg is to create the situation for AT with the maximal possible entropy.\par
\vspace{.2cm}
\section{The Full Solution of the Game $G(n=2,k=1, m=1, a,b)$}
As an example of how Theorem 3 works, we outline the full solution of game $G(n=2,k=1, m=1)$ with general vector of
values $c=(c, 1), c\geq 1$ and symmetric sensitivities/specificities $a_1=a_2=a, b_1=b_2=b, 1/2< a,b <1$. WLOG the \prb of explosion  $p$ can be assumed to equal $1$. In this game, there are two possible positions of locks: $\gamma_1=(0,1)$
and $\gamma_2=(1,0)$, with prior \dst $b(\gamma_1)=x, 0\leq x \leq 1$, and four signals $s=(s_1,s_2)\equiv \nu$:
$\nu_1=(-,+)\equiv (0,1), \nu_2=(+,-), \nu_3=(-,-), \nu_4=(+,+)$. The conditional \prbs of these signals
$b(s|\gamma_1 )\equiv b(\nu_i|\gamma_1), i=1,2,3,4$ are $ba=e_1, (1-b)(1-a)=e_2, b(1-a)=e_3, (1-b)a=e_4$, and the conditional \prbs $b(\nu_i|\gamma_2 )$ are $e_2,e_1,e_3,e_4$. Under our assumptions about $a,b$, we have $e_2<e_3,e_4< e_1$. To find the best response of AT for each value $x$ and each signal $s\equiv \nu_i$, we have to compare the \emph{potential} damage $d_1(s|x)$ to site 1 if a bomb is placed there, equal to $c\ast P(T_1=0|s)$, with the similar potential damage to site 2, $d_2(s|x)$, equal to $1\ast P(T_2=0|s)$. The \emph{real} damage $d(s|x)$ with optimal placement of a unique bomb is given by the maximum of these two functions for each $s$ and $x$. We denote $P(T_i=0|s)=\alpha_i(s)$. \par
We have $\alpha_1(s)=b(\gamma_1|s)=b(\gamma_1)b(s|\gamma_1)/p(s)$ and $\alpha_2(s)=
b(\gamma_2|s)=b(\gamma_2)b(s|\gamma_2)/p(s).$ Since $b(\gamma_1)=x,b(\gamma_2)=1-x$, for signal $\nu_1$ we have $\alpha_1(\nu_1)=xe_1/p(\nu_1)$ and $\alpha_2(\nu_1)=(1-x)e_2/p(\nu_1)$. Thus, the decision as to where to place the bomb is determined by comparing  $d_1(\nu_1|x)=cxe_1/p(\nu_1)$ with $d_2(\nu_1|x)=(1-x)e_2/p(\nu_1)$. Let us denote $\rho_1$ the value of $x$ where $d_i(\nu_1|x)$ are equal. Solving the equality $cxe_1=(1-x)e_2$, we obtain that $\rho_1=1/(ce_1/e_2+1).$  This means that the real damage $d(\nu_1|x)=(1-x)e_2/p(\nu_1)$ on the interval $(0, \rho_1)$, and $d(\nu_1|x)=cxe_1/p(\nu_1)$ on the interval $(\rho_1, 1)$. \par
Similarly, for signal $\nu_2$, we have to compare function $cxe_2$ with function $(1-x)e_1$. The second function is larger than the first one on the interval $(0, \rho_2)$, where $\rho_2=1/(ce_2/e_1+1).$ 
For both signals $\nu_3,\nu_4$ we have to compare function $cx$ with function $1-x$. The second function is larger than the first one on the interval $(0, \rho_3)$, where $\rho_3=1/(c+1).$ Since $e_1/e_2>1$, we have $\rho_1<\rho_3<\rho_2<1.$\par
As a result, the interval $(0,1)$ is partitioned into four disjoint subintervals $\Delta_1=(0, \rho_1), \Delta_2=(\rho_1, \rho_3),\Delta_3=(\rho_3, \rho_2), \Delta_4=(\rho_2, 1)$.
On subinterval $\Delta_1$, where with high \prb $1-x$ the lock is in site 1, for all signals the bomb is placed into site 2; on $\Delta_4$, where $x$ is closer to 1, the situation is reversed, and for all signals the bomb is placed into site 1; on $\Delta_2$ for signal $\nu_1$ the bomb moves to site 1; and for the three other signals it remains in site 2;
on $\Delta_3$ for signal $\nu_2$, the bomb remains in site 2; and for the three other signals it moves to site 1. As a result, the function of the real damage $d(x)=\sum_{s}p(s|x)d(s|x)$ is $(1-x)$ on interval $\Delta_1$, is $cx$ on interval $\Delta_4$, and on the intervals $\Delta_2, \Delta_3$ has the \flw form:
\begin{eqnarray}
d(x)&=&e_1[\frac{e_1+e_3+e_4}{e_1}+x(c-\frac{e_1+e_3+e_4}{e_1})]\ \ \text{on interval}\  \Delta_2, \notag \\
    &=&e_1[1+x(c\frac{e_1+e_3+e_4}{e_1}-1)] \ \ \text{on interval}\  \Delta_3.  \label{dx}
\end{eqnarray}
Let  $c_*=\ds \frac{e_1+e_3+e_4}{e_1}\equiv \frac{a+b}{ab}-1>1$. Formula (\ref{dx}) implies that function $d(x)$ is decreasing on interval $\Delta_1$, increasing on intervals $\Delta_3, \Delta_4$, and if \ $1\leq c<c_*$, then $d(x)$ is decreasing on interval $\Delta_2$, and hence the minimum of function $d(x)$ is reached at point $x_*=\rho_3=1/(c+1)$ and $d(x_*)=d(\rho_3)$. If $c>c_*$, then $d(x)$ is increasing on interval $\Delta_2$, and hence the minimum of function $d(x)$ is reached at point $x_*=\rho_1=e_2/(ce_1+e_2)$ and $d(x_*)=d(\rho_1)$. If $c=c_*$, then the minimum  of function $d(x)$ is reached on the whole interval $\Delta_2$  and $d(x_*)=d(\rho_1)=d(\rho_3)$. Thus, we obtained the
best \stg of DF, $b_*(\gamma_1)\equiv x_*=x_*(c, e_1,e_2,e_3,e_4)$ and the value of the game $v_*=\sum_{s}p(s|x_*)d(s|x_*)$. They both can be expressed in explicit form, as a function of  parameters $c, e_1,e_2,e_3,e_4$ and $k=e_1/e_2\equiv ab/(1-a)(1-b)$, as \fls:
\vspace{-.2cm}
\begin{eqnarray}
x_*&=&\rho_3=\frac{1}{1+c}, \ \  v_*=d(\rho_3)=\frac{c(1+e_1-e_2)}{1+c}\equiv \frac{c(a+b)}{1+c} \  \ \text{if  }\ \ c\leq c_*, \notag \\
&=&\rho_1=\frac{1}{1+ck}, \ \  v_*=d(\rho_1)=\frac{c}{e_2/e_1+c}=\frac{c}{1/k+c} \ \ \ \text{if  }\ \ c\geq c_*.  \label{v*}
\end{eqnarray}
For $c=c_*$ we obtain $v_*=d(\rho_1)=d(\rho_3)=a+b-ab$.\\
The best \stg of AT depends on signal $s$, and is defined by functions $d_i(s|x_*)$.\par
Note that when the testing is perfect, $a_i=b_i=1$, we have $e_1=1, e_2=e_3=e_4=0$ and obtain that $x_*=\rho_1=0, v_*=1$. In other words, the only lock is placed into more valuable site 1, and the loss is equal to 1. If the testing is non-informative, $a_i=b_i=1=1/2$, we have $e_1=e_2=e_3=e_4=1/4$ and obtain that $x_*=\frac{1}{1+c}$ and $v_*=\frac{c}{1+c}.$\par
For the values $a=7/12, b=9/12$, we obtain $e_i=(63, 15, 45, 21)/144, k=\frac{21}{5}=4.2, c_*=129/63\approx 2.048$. Then e.g. for $c=2$ we have $x_*=\frac{1}{3}, v_*=\frac{8}{9}\approx .889$. For for $c=3$ we have $x_*=\frac{1}{1+3*21/5}=\frac{5}{68}\approx .074, v_*=\frac{63}{68}\approx .926$.\par
It is interesting to note that in this example for $1\leq c\leq c_*$ the optimal \stg of DF is the same as in the case of noninformative case, i.e. $\alpha =(\frac{1}{1+c}, \frac{c}{1+c})$ but this equality is not true for $c>c_*$. This tells that the Theorem 3 in this example is working, i.e. for optimal \stg $b_*(\gamma)\equiv x_*$ the equality $d_{1,*}=d_{1,*}$, the Indifference Principle holds, but the equality $c_i\alpha_i=const$ does not. \par
This relatively simple example shows also that the exact solution of a general $G$ game can be a challenging problem.\par

}

\end{document}